\documentclass [11pt] {article}

\usepackage{amsfonts}
\usepackage{amsmath}
\usepackage{amssymb}
\usepackage{color}
\usepackage{graphicx}
\usepackage{times}
\usepackage[square, numbers]{natbib}
\usepackage[T1]{fontenc}
\usepackage{url}
\usepackage{setspace}

\newcommand {\cse} {CSE}

\begin{document}

\title{Computational Science and Engineering}
\title{\hspace*{-1em}Computational Science, and Engineering, Education\hspace*{-1em}}
\title{Computational Science Education}
\title{\hspace*{-1em}Education in Computational Science and Engineering\hspace*{-1em}}
\title{\hspace*{-1em}Education for Computational Science and Engineering\hspace*{-1em}}

\author{Joseph F. Grcar\,\footnote{6059 Castlebrook Drive; Castro Valley, CA 94552-1645 USA} \,\footnote {\texttt {jfgrcar@comcast.net}, \texttt {jfgrcar@gmail.com}.}}

\date{}

\maketitle

\onehalfspace

\begin {abstract}
Computational science and engineering (CSE) has been misunderstood to advance with the construction of enormous computers. To the contrary, the historical record demonstrates that innovations in CSE come from improvements to the mathematics embodied by computer programs. Whether scientists and engineers become inventors who make these breakthroughs depends on circumstances and the interdisciplinary extent of their educations. The USA currently has the largest CSE professorate, but the data suggest this prominence is ephemeral. 
 \end {abstract}

Computational science and engineering (CSE) embodies the challenges in President Obama's vision of the future, ``ours to win'' \citep {POTUS2011}. Wealthy nations with comparable scientific institutions can achieve parity in technical prowess. When countries have have similar populations of technical professionals, one way to excel is to make better use of expert talent through computing. Research and development (R\&D) proceeds more quickly and with greater certainty by using computers to analyze information and to pursue the implications of scientific theories. 
 
The preceding ``elevator speech'' has been heard in Washington before.\footnote {For a brief history of CSE, search the text of articles in \textit {Nature\/} or \textit {Science\/} for the exact phrase ``computational science.'' For a more detailed timeline, search for ``supercomputer.''} From ``NSF Plans Help with Big Computer Problems'' \citep {Walsh1984} in 1985, through ``Energy Labs Urged to Boost Supercomputing Capability'' \citep {Macilwain1997} in 1997, to ``An Endless Frontier Postponed'' \citep {Lazowska2005} in 2005, \cse\ has been used to justify the purchase of massive computers --- and is being used again, in President Obama's budget. To the contrary, the fastest computers are irrelevant to the innovations in \cse\ that bettered peoples' lives and transformed how R\&D is done. 

The breakthroughs in \cse\ are the mathematical programs of computing, not the machines themselves. President Obama misunderstood the ``Sputnik moment'' when the People's Republic of China (PRC) built ``the world's fastest computer.'' Japan and the USA previously built such machines to little effect because availability of the fastest computers is severely limited by cost \citep {Sun1989}. The cost of the absolutely fastest machines remains high even though their composition has changed. The extreme engineering that once made computers fast was overtaken by silicon miniaturization for the mass market. Today's fastest machines are merely huge ``clusters'' of commodity devices working in parallel \citep {Hillis1993, Khamsi2004}, which many research groups can acquire in small configurations \citep {Hargrove2001}. These cluster computers were predicted to be difficult to use \citep [p.\ 458] {Barker1994}, and have proven so. 

At present there are two paths to achieving great speed. Either computer programs might be written to excel on clusters, or special chips might be built for one type of calculation. (Examples of specialty chips are in cell phones and game consoles.) Programming skill is also critical for the special-purpose chips because likely programs must be well characterized before building silicon to suit. The scientific literature has examples of both approaches for calculations of molecular dynamics \citep {Bowers2006, Shaw2007, Service2010}. The flexing motions of large molecules affect their reactivity \citep {Ammal2003}, so these calculations reveal how molecules are biologically active. 

The importance of clever programming begins only after there is something to write a program about. The prerequisite mathematical invention often is obscured by the continuing work to keep up with computing practice. For example, at the turn of the 19th century, A. M. Legendre and C.\ F.\ Gauss \citep {Legendre1805-HM, Gauss1809-HM} invented what is now called regression analysis, which is used to infer parametric models in diverse fields such as econometrics and epidemiology. Gauss in 1810 \citep {Gauss1810-HM} also invented a way to calculate regression coefficients by hand that was used for 150 years, with modifications for various types of manual computing \citep {Grcar2011c, Grcar2011e}. G.\ H.\ Golub in 1965 \citep {Golub1965} invented a better way for electronic computers. Today, Golub's method is found in libraries of computer programs, where it evolves with programming techniques \citep {Nerlove2004, Trefethen2007}.

The example of regression analysis illustrates that the paramount breakthroughs in \cse\ are the mathematical inventions that live through generations of computer programs. These advances come from researchers who learn to work with mathematics across disciplinary boundaries \citep {Grcar2011g}. Prominent examples are: (a) a fast method to evaluate convolutions was invented by two polymaths working in statistics and computer science  \citep {Brillinger2002}; (b) a method for inferring structures of proteins from x-ray diffraction patterns was invented jointly by a physicist and a mathematician \citep {Gilmore1980, Karle1986, Hauptman1986}; and (c) a method for identifying the current state of a dynamical system was invented by a control engineer and refined by a programmer for the tiny computers on early space vehicles \citep {Grewal2010}. Some of these inventions enable scientific research, while others are ubiquitous in daily life; the importance of the inventions (a, b, c) is discussed in \citep {Bracewell1990, Hendrickson1986, NAE2008}, respectively. 

Whether scientists and engineers become inventors who make \cse\ breakthroughs depends on circumstances and their educations. Success in President Obama's challenge is not decided by the fastest computer, but it may be decided by the size and quality of the CSE professorate. To that end, a census was made of educators who are experts about calculations and the mathematics of a science or engineering field. (The census is described in the Appendix.)

\begin {table} [b]
\caption {\it (Column 1) populations of educators specializing in \cse, (column 2) general populations, (column 3) ratio of \cse\ to general. Countries are ordered by decreasing ratios. \textbf {Bold} are countries with more CSE educators per-capita than the USA. Many of these countries (*) are members of the EU; Hong Kong is an administrative region of the PRC. Ratios for Japan and PRC may be suppressed by undercounting the scientific literature written in non-Latin alphabets.}
\label {tab:countries}
\begin {center}
\footnotesize
\begin {tabular} {| r | r r c |}
\cline {2-4}
\multicolumn {1} {c} {\vrule depth0pt height2.0ex width 0pt}& \multicolumn {2} {| c} {\textsc {\ populations}}&\\
\multicolumn {1} {c} {\vrule depth0pt height2.0ex width 0pt}& \multicolumn {2} {| c} {\textsc {(\,usa = 100\,)}}&\\
\multicolumn {1} {c} {}& \multicolumn {1} {| r} {\textsc {cse}}& \multicolumn {1} {c} {\textsc {gen.}}&\multicolumn {1} {c |} {\textsc {ratio}}\\ 
\multicolumn {1} {c} {}& \multicolumn {1} {| r} {\textsc {educ.}}& \multicolumn {1} {c} {\textsc {pop.}}& \multicolumn {1} {c |} {\textsc {cse/gen}}\\ \hline
\vrule depth0pt height2.5ex width 0pt
\textbf {Hong Kong}& 12.3& 2.3& 5.4\\
\textbf {Israel}& 13.4& 2.5& 5.4\\
\textbf {Singapore}& 5.6& 1.6& 3.5\\
$^*$\textbf {Austria}& 8.4& 2.7& 3.1\\
\textbf {Switzerland}& 5.0& 2.5& 2.0\\ 
\textbf {Canada}& 19.0& 11.1& 1.7\\
$^*$\textbf {Germany}& 40.8& 26.5& 1.5\\
$^*$\textbf {France}& 31.3& 21.3& 1.5\\
$^*$\textbf {United Kingdom}& 26.8& 20.1& 1.3\\ 
\vrule depth0pt height3ex width0pt USA& 100\phantom{.0}& 100\phantom{.0}& \ 1\phantom{\,\,.0}\\
Australia& 7.3& 7.3& 1.0\\
EU combined& 165.4& 162.3& 1.0\\
$^*$Netherlands& 5.6& 5.4& 1.0\\
$^*$Italy& 15.6& 19.6& 0.8\\
$^*$Spain& 10.1& 14.9& 0.7\\
Japan& 12.3& 41.3& 0.3\\
PRC& 54.2& 433.6& 0.1\\
\hline
\end {tabular}
\end {center}
\end {table}

The USA currently has the largest \cse\ professorate (column 1 in \textbf {Table \ref {tab:countries}}), but the data suggest this prominence is likely to decline. The European Union (EU) already has many more \cse\ educators than the USA and the same number per capita (column 3). European and Asian countries known for exporting manufactured goods have more \cse\ educators per-capita than the USA; in contrast, the USA exports mostly agricultural commodities. While the People's Republic of China (PRC) continues to develop a manufacturing economy, it already has the second largest number of \cse\ educators of any country.

\begin {figure} [h]
\begin {center}
\smallskip
\includegraphics [scale=0.8] {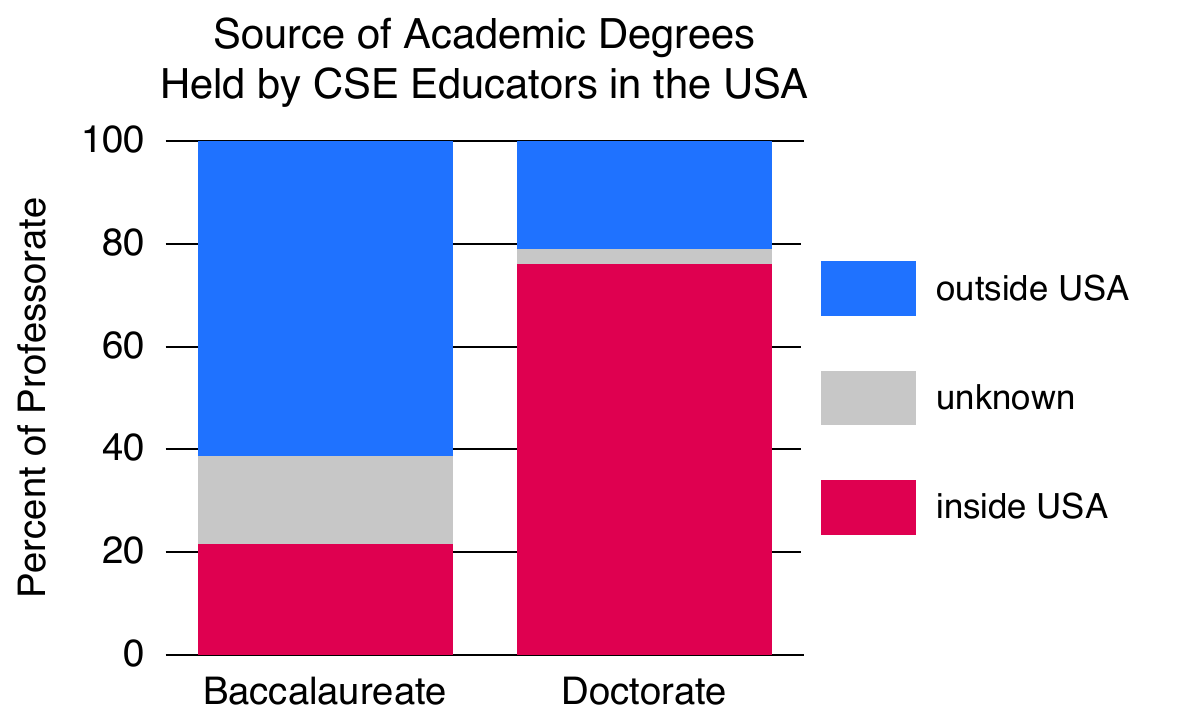}
\end {center}
\vspace*{-3ex}
\caption {\it Country of origin for academic degrees held by \cse\ educators in the USA. Most \cse\ educators earned baccalaureate degrees outside the country, which suggests that the immigration of academically talented individuals is critical to the continued prowess of the USA in CSE.}
\label {fig:source}
\end {figure}

The greatest influence on \cse\ education is the supply of educators.  Whether the USA can sustain the present level of effort is moot because over half the current professorate earned the baccalaureate degree outside the country (\textbf {Figure \ref {fig:source}}). As enrollments by foreign students fluctuate \citep {Bhattacharjee2005}, preeminence in CSE for the USA may depend on extraneous factors such as immigration policy.

What universities do control is the content and quality of education. \cse\ poses a unique combination of difficulties. The first is, knowledge of using computers in a scientific or engineering field does not suffice to extend the methodology. As seen in the prominent examples above, advances in \cse\ require a combination of mathematics, computer science, and disciplinary (or domain-specific) training. Second is,  academic departments for mathematics or computer science do not serve as nuclei for \cse\ education. Elite mathematics departments in the USA shun interdisciplinary subjects \citep {Davis1994-complete, Grcar2010e, Grcar2011b}. Computer science departments do have interdisciplinary research interests, but mostly apart from \cse\ \citep {Shaw2004, Kowalik2006}. Mathematics and computer science departments also have declining matriculations \citep {Foster2005, Bressoud2009}, which impedes accepting new responsibilities. Third is,  interdisciplinary educational programs do not automatically follow from interdisciplinary research programs. University administrators in the USA are adept at encouraging interdisciplinary research \citep {Sa2008}. For example, the National Science Foundation built computer centers to support research faculty in the Cold War. However, early in the age of cheap computers it was foreseen that \cse\ education would have to be addressed directly \citep {Fox1992, Sameh1993}. 

\begin {figure} [h]
\begin {center}
\smallskip
\includegraphics [scale=0.8] {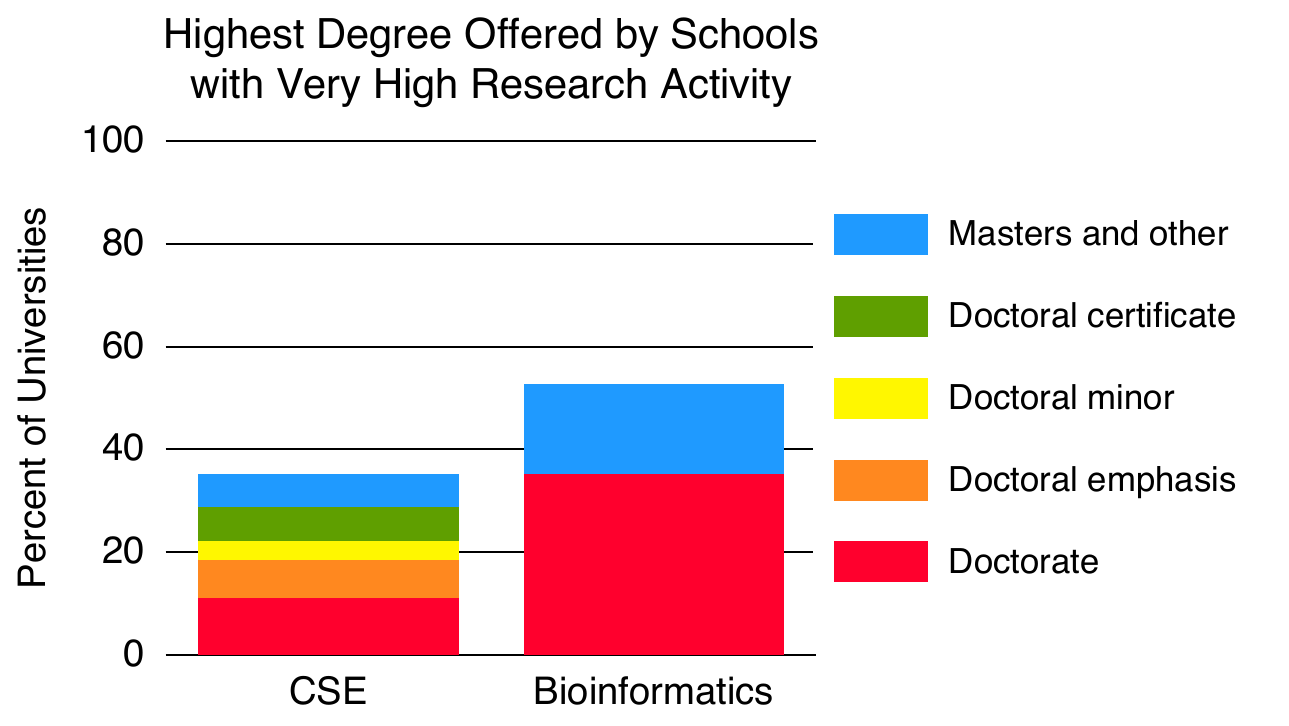}
\end {center}
\vspace*{-3ex}
\caption {\it Highest degrees in \cse\ and in bioinformatics offered by the 108 universities in the USA with very high research activity. Many more schools have established education programs for bioinformatics than for \cse.}
\label {fig:programs}
\end {figure}

As a consequence of the hurdles to establishing educational programs for \cse, the 108 universities in the USA with very high research activity \citep {Carnegie2010} have a total of just 37 programs for \cse\ education (\textbf {Figure \ref {fig:programs}}). Illustrating the diversity of administrative structures: a few programs are departments, e.g.\ \citep {GeorgiaTech}, while others are interdepartmental institutes with their own degree programs, e.g.\ \citep {Stanford}, and still others are subjects of emphasis for degrees in participating departments, e.g.\ \citep {Purdue, SantaBarbara}. These programs are discussed in a small and thinly cited literature \citep {SIAM2001, Yasar2003, Manson2010}. In comparison, the field of bioinformatics is younger than CSE and is already well represented by degree programs and departments in schools of medicine and public health \citep {Watanabe2004}. Considering the historical contributions that \cse\ has made (which enable all manner of commercial, consumer, medical, military, and scientific devices), the future appears to be the United States' to lose.

\appendix

\section* {Appendix}

This appendix describes how data were gathered for \textbf {Table \ref {tab:countries}} and \textbf {Figure \ref {fig:source}}.

Educators are identified from their publications in the database of \textit {Zentralblatt f\"ur Mathematik und ihre Grenzgebiete\/} \citep {Zentralblatt-fur-Mathematik, Fairweather2009}. The selection criteria are that publications belong to topical classifications 65 or 68 (numerical analysis or computer science) and simultaneously one or more of 60, 62, 70--94 (applied mathematics of biology, control, information, mechanics, optics, statistics, etc.). Most publications appear outside core mathematics journals as indicated by the database name. Educators are chosen with ten papers in the decade 2000--2009 as evidence of recent, continuing research and to keep the sample manageable. These criteria identify 766 individuals, all of whom are associated with academic institutions.

For perspective on the sample size of 766 CSE educators with ten research publications: (a) the Society for Industrial and Applied Mathematics has a \cse\ activity group with 936 non-student members, (b) the entire field of probability and statistics (which underlies much research methodology in life sciences) has 2826 authors with ten or more publications, and (c) \cse\ has 3546 authors with five or more publications.

{
\raggedright

\begin{thebibliography}{52}
\providecommand{\natexlab}[1]{#1}
\providecommand{\url}[1]{\texttt{#1}}
\expandafter\ifx\csname urlstyle\endcsname\relax
  \providecommand{\doi}[1]{doi: #1}\else
  \providecommand{\doi}{doi: \begingroup \urlstyle{rm}\Url}\fi

\bibitem[Obama(2011)]{POTUS2011}
B.~H. Obama.
\newblock Remarks by the {P}resident in {S}tate of {U}nion {A}ddress.
\newblock White House press release, January 25, 2011.

\bibitem[Walsh(1984)]{Walsh1984}
J.~Walsh.
\newblock {NSF} plans help with big computer problems.
\newblock \emph{Science}, 223\penalty0 (4638):\penalty0 797--798, 1984.
\newblock \doi{10.1126/science.223.4638.797}.

\bibitem[Macilwain(1997)]{Macilwain1997}
C.~Macilwain.
\newblock Energy labs urged to boost supercomputing capability.
\newblock \emph{Nature}, 390\penalty0 (6661):\penalty0 651, 1997.
\newblock \doi{10.1038/37698}.

\bibitem[Lazowska and Patterson(2005)]{Lazowska2005}
E.~D. Lazowska and D.~A. Patterson.
\newblock An endless frontier postponed.
\newblock \emph{Science}, 308\penalty0 (5723):\penalty0 757, 2005.
\newblock \doi{10.1126/science.1113963}.

\bibitem[Sun(1989)]{Sun1989}
M.~Sun.
\newblock Supercomputer market needs supersalesmen.
\newblock \emph{Science}, 245\penalty0 (4918):\penalty0 596--597, 1989.
\newblock \doi{10.1126/science.245.4918.596}.

\bibitem[Hillis and Boghosian(1993)]{Hillis1993}
W.~D. Hillis and B.~M. Boghosian.
\newblock Parallel scientific computation.
\newblock \emph{Science}, 261\penalty0 (5123):\penalty0 856--863, 1993.
\newblock \doi{10.1126/science.261.5123.856}.

\bibitem[Khamsi(2004)]{Khamsi2004}
R.~Khamsi.
\newblock A speedy recipe for supercomputing.
\newblock \emph{Nature News}, 2004.
\newblock \doi{10.1038/news041108-3}.

\bibitem[Hargrove et~al.(2001)Hargrove, Hoffman, and Sterling]{Hargrove2001}
W.~W. Hargrove, F.~M. Hoffman, and T.~Sterling.
\newblock The do-it-yourself supercomputer.
\newblock \emph{Scientific American}, 285:\penalty0 72--79, 2001.
\newblock \doi{10.1038/scientificamerican0801-72}.

\bibitem[Barker(1994)]{Barker1994}
B.~Barker.
\newblock Lessons learned.
\newblock In K.~R. Ames and A.~Brenner, editors, \emph{Frontiers of
  Supercomputing II: A National Reassessment}, pages 453--468. Univ.\ of
  California Press, 1994.

\bibitem[Bowers and coauthors(2006)]{Bowers2006}
K.~J. Bowers and 12~coauthors.
\newblock Scalable algorithms for molecular dynamics simulations on commodity
  clusters.
\newblock In \emph{Proc.\ 2006 ACM/IEEE Conf.\ on Supercomputing}. 2006.
\newblock \doi{10.1145/1188455.1188544}.

\bibitem[Shaw and coauthors(2007)]{Shaw2007}
D.~E. Shaw and 26~coauthors.
\newblock Anton, a special-purpose machine for molecular dynamics simulation.
\newblock \emph{ACM SIGARCH Computer Architecture News}, 35\penalty0
  (2):\penalty0 1--12, 2007.
\newblock \doi{10.1145/1273440.1250664}.

\bibitem[Service(2010)]{Service2010}
R.~F. Service.
\newblock Custom-built supercomputer brings protein folding into view.
\newblock \emph{Science}, 330\penalty0 (6002):\penalty0 308--309, 2010.
\newblock \doi{10.1126/science.330.6002.308-a}.

\bibitem[Ammal and coauthors(2003)]{Ammal2003}
S.~C. Ammal and 3~coauthors.
\newblock Dynamics-driven reaction pathway in an intramolecular rearrangement.
\newblock \emph{Science}, 299\penalty0 (5612):\penalty0 1555--1557, 2003.

\bibitem[Legendre(1805)]{Legendre1805-HM}
A.~M. Legendre.
\newblock \emph{Nouvelle m\'ethodes pour la d\'etermination des orbites des
  com\`etes}.
\newblock Chez Didot, Paris, 1805.

\bibitem[Gauss(1809)]{Gauss1809-HM}
C.~F. Gauss.
\newblock \emph{Theoria Motus Corporum Coelestium in Sectionibus Conicis Solum
  Ambientium}.
\newblock Perthes and Besser, Hamburg, 1809.

\bibitem[Gauss(1810)]{Gauss1810-HM}
C.~F. Gauss.
\newblock Disquisitio de elementis ellipticis {P}alladis.
\newblock \emph{Commentationes Societatis Regiae Scientiarum Gottingensis
  recentiores: Commentationes classis mathematicae}, 1 (1808--1811):\penalty0
  1--26, 1810.

\bibitem[Grcar(2011{\natexlab{a}})]{Grcar2011c}
J.~F. Grcar.
\newblock How ordinary elimination became {G}aussian elimination.
\newblock \emph{Historia Math.}, 38\penalty0 (2):\penalty0 163--218,
  2011{\natexlab{a}}.
\newblock \doi{10.1016/j.hm.2010.06.003}.

\bibitem[Grcar(2011{\natexlab{b}})]{Grcar2011e}
J.~F. Grcar.
\newblock Mathematicians of {G}aussian ellimination.
\newblock \emph{Notices Amer. Math. Soc.}, 58\penalty0 (6):\penalty0 782--792,
  2011{\natexlab{b}}.

\bibitem[Golub(1965)]{Golub1965}
G.~H. Golub.
\newblock Numerical methods for solving linear least squares problems.
\newblock \emph{Numer. Math.}, 7:\penalty0 206--216, 1965.

\bibitem[Nerlove(2004)]{Nerlove2004}
M.~Nerlove.
\newblock Programming languages: A short history for economists.
\newblock \emph{J.\ Econ.\ Soc.\ Meas.}, 29\penalty0 (1--3):\penalty0 189--203,
  2004.

\bibitem[Trefethen(2007)]{Trefethen2007}
L.~N. Trefethen.
\newblock {Obituary: Gene H. Golub (1932--2007)}.
\newblock \emph{Nature}, 450\penalty0 (7172):\penalty0 962, 2007.
\newblock \doi{10.1038/450962a}.

\bibitem[Grcar(2011{\natexlab{c}})]{Grcar2011g}
J.~F. Grcar.
\newblock Great ideas of mathematics of value to society.
\newblock submitted, 2011{\natexlab{c}}.

\bibitem[Brillinger(2002)]{Brillinger2002}
D.~R. Brillinger.
\newblock {John Wilder Tukey (1915--2000)}.
\newblock \emph{Notices Amer. Math. Soc.}, 49\penalty0 (2):\penalty0 193--201,
  2002.

\bibitem[Gilmore(286)]{Gilmore1980}
C.~J. Gilmore.
\newblock Direct methods and the phase problem.
\newblock \emph{Nature}, pages 315--316, 286.
\newblock \doi{10.1038/286315a0}.

\bibitem[Karle(1986)]{Karle1986}
J.~Karle.
\newblock Recovering phase information from intensity data.
\newblock \emph{Science}, 232\penalty0 (4752):\penalty0 837--843, 1986.
\newblock \doi{10.1126/science.232.4752.837}.

\bibitem[Hauptman(1986)]{Hauptman1986}
H.~Hauptman.
\newblock The direct methods of x-ray crystallography.
\newblock \emph{Science}, 233\penalty0 (4760):\penalty0 178--183, 1986.
\newblock \doi{10.1126/science.233.4760.178}.

\bibitem[Grewal and Andrews(2010)]{Grewal2010}
M.~S. Grewal and A.~P. Andrews.
\newblock Applications of {K}alman filtering in aerospace 1960 to the present.
\newblock \emph{IEEE Contr. Syst, Mag.}, 3:\penalty0 69--78, 2010.

\bibitem[Bracewell(1990)]{Bracewell1990}
R.~N. Bracewell.
\newblock Numerical transforms.
\newblock \emph{Science}, 248\penalty0 (4956):\penalty0 697--704, 1990.
\newblock \doi{10.1126/science.248.4956.697}.

\bibitem[Hendrickson(1986)]{Hendrickson1986}
W.~A. Hendrickson.
\newblock The 1985 {N}obel prize in chemistry.
\newblock \emph{Science}, 231\penalty0 (4736):\penalty0 362--364, 1986.
\newblock \doi{10.1126/science.231.4736.362}.

\bibitem[{Nat.\ Acad.\ Eng.\ USA}(2008)]{NAE2008}
{Nat.\ Acad.\ Eng.\ USA}.
\newblock {Charles Stark Draper Prize}.
\newblock Press Release, 2008.

\bibitem[Bhattacharjee(2005)]{Bhattacharjee2005}
Y.~Bhattacharjee.
\newblock Schools cheer rise in foreign students.
\newblock \emph{Science}, 310\penalty0 (5750):\penalty0 957, 2005.

\bibitem[Davis(1994)]{Davis1994-complete}
C.~Davis.
\newblock Where did twentieth-century mathematics go wrong?
\newblock In \emph{The Intersection of History and Mathematics}, pages
  129--142. Birkh\"auser, Basel, 1994.

\bibitem[Grcar(2010)]{Grcar2010e}
J.~F. Grcar.
\newblock Topical bias in generalist mathematics jounals.
\newblock \emph{Notices Amer. Math. Soc.}, 57\penalty0 (11):\penalty0
  1421--1424, 2010.

\bibitem[Grcar(2011{\natexlab{d}})]{Grcar2011b}
J.~F. Grcar.
\newblock Mathematics turned inside out: the intensive faculty versus the
  extensive faculty.
\newblock \emph{Higher Education}, 61\penalty0 (6):\penalty0 693--720,
  2011{\natexlab{d}}.
\newblock \doi{10.1007/s10734-010-9358-y}.

\bibitem[Shaw~et al.(2004)]{Shaw2004}
M.~Shaw~et al.
\newblock \emph{Computer Science: Reflections on the Field, Reflections from
  the Field}.
\newblock The National Academies Press, Washington, 2004.

\bibitem[Kowalik(2006)]{Kowalik2006}
J.~Kowalik.
\newblock The applied mathematics and computer science schism.
\newblock \emph{Computer}, 39\penalty0 (3):\penalty0 102--104, 2006.
\newblock \doi{10.1109/MC.2006.106}.

\bibitem[Foster(2005)]{Foster2005}
A.~L. Foster.
\newblock Student interest in computer science plummets.
\newblock \emph{The Chronicle of Higher Education}, 51\penalty0 (38):\penalty0
  A31, May 27 2005.

\bibitem[Bressoud(2009)]{Bressoud2009}
D.~M. Bressoud.
\newblock Is the sky still falling?
\newblock \emph{Notices Amer. Math. Soc.}, 56\penalty0 (1):\penalty0 20--25,
  2009.

\bibitem[S\'{a}(2008)]{Sa2008}
C.~M. S\'{a}.
\newblock `{I}nterdisciplinary strategies' in {U}.{S}.\ research universities.
\newblock \emph{Higher Education}, 55:\penalty0 537--552, 2008.
\newblock \doi{10.1007/s10734-007-9073-5}.

\bibitem[Fox(1992)]{Fox1992}
G.~C. Fox.
\newblock Parallel computing and education.
\newblock \emph{Daedalus}, 121\penalty0 (1):\penalty0 111--118, 1992.

\bibitem[Sameh and Riganati(1993)]{Sameh1993}
A.~H. Sameh and J.~Riganati.
\newblock Computational science and engineering.
\newblock \emph{IEEE Computer}, 26\penalty0 (10):\penalty0 8--12, 1993.
\newblock \doi{10.1109/MC.1993.10105}.

\bibitem[{Carnegie Foundation}(2010)]{Carnegie2010}
{Carnegie Foundation}.
\newblock The {C}arnegie classification of institutions of higher
  education\texttrademark, 2010.
\newblock URL \url{http://classifications.carnegiefoundation.org/}.

\bibitem[{Georgia Institute of Technology}()]{GeorgiaTech}
{Georgia Institute of Technology}.
\newblock {School of Computational Science and Engineering}.
\newblock URL \url{http://www.cse.gatech.edu/}.

\bibitem[{Stanford University}()]{Stanford}
{Stanford University}.
\newblock {Institute for Computational and Mathematical Engineering}.
\newblock URL \url{http://icme.stanford.edu/}.

\bibitem[{Purdue University}()]{Purdue}
{Purdue University}.
\newblock {Computational Science and Engineering Graduate Program}.
\newblock URL \url{http://www.gradschool.purdue.edu/cse/}.

\bibitem[{University of California-Santa Barbara}()]{SantaBarbara}
{University of California-Santa Barbara}.
\newblock {Computational Science and Engineering}.
\newblock URL \url{http://www.cse.ucsb.edu/}.

\bibitem[{SIAM Working Group on CSE Education}(2001)]{SIAM2001}
{SIAM Working Group on CSE Education}.
\newblock Graduate education in computational science and engineering.
\newblock \emph{SIAM Rev.}, 43\penalty0 (1):\penalty0 163--177, 2001.
\newblock \doi{10.1137/S0036144500379745}.

\bibitem[Ya\c{s}ar and Landau(2003)]{Yasar2003}
O.~Ya\c{s}ar and R.~H. Landau.
\newblock Elements of computational science and engineering education.
\newblock \emph{SIAM Rev.}, 45\penalty0 (4):\penalty0 787--805, 2003.
\newblock \doi{10.1137/S0036144502408075}.

\bibitem[Manson and Olsen(2010)]{Manson2010}
J.~R. Manson and R.~J. Olsen.
\newblock Assessing and refining an undergraduate computational science
  curriculum.
\newblock \emph{Procedia Computer Science}, 1\penalty0 (1):\penalty0 857--865,
  2010.
\newblock \doi{10.1016/j.procs.2010.04.094}.

\bibitem[Watanabe(2004)]{Watanabe2004}
M.~Watanabe.
\newblock Changing of the guard.
\newblock \emph{Nature}, 428\penalty0 (6982):\penalty0 584--585, 2004.
\newblock \doi{10.1038/nj6982-584a}.

\bibitem[Zen()]{Zentralblatt-fur-Mathematik}
\emph{Zentralblatt f\"ur Mathematik}.
\newblock Springer-Verlag, Berlin.
\newblock \url {http://www.zentralblatt-math.org/zmath/en/}.

\bibitem[Fairweather and Wegner(2009)]{Fairweather2009}
G.~Fairweather and B.~Wegner.
\newblock {Mathematics Subject Classification 2010}.
\newblock \emph{Notices Amer. Math. Soc.}, 56\penalty0 (7):\penalty0 848, 2009.

\end{thebibliography}

}

\end {document}